\newif\ifsiam
\siamfalse
\ifdefined\SIAMVERSION\siamtrue\fi

\ifsiam
    \documentclass[review,onefignum,onetabnum,onealgnum]{siamart250211}
\else
    \documentclass[11pt]{article}
    \usepackage[margin=1in]{geometry}
    \usepackage{algorithm}
\fi
\usepackage{algpseudocode}
\usepackage{url}

\input{splitted/preamble}
\input{splitted/shared}

\title{Mixed-Precision Coalesced Matrix-Free Multigrid with Vertex-Patch Smoothing}

\ifsiam
    \headers{Mixed-Precision Vertex-Patch Multigrid}{M. Wichrowski}
    \author{Michał Wichrowski\thanks{\email{mwichro@mimuw.edu.pl}}}
\else
    \author{Michał Wichrowski}
    \date{}
    \newenvironment{keywords}{\par\noindent\textbf{Keywords:}}{\par}
    \newenvironment{MSCcodes}{\par\noindent\textbf{AMS subject classifications:}}{\par}
\fi
\renewcommand{\paragraph}[1]{\par\textbf{#1.}\ \ignorespaces}
\begin{document}

\maketitle
\ifsiam\else
    \footnotetext{\texttt{mwichro@mimuw.edu.pl}}
\fi

\begin{abstract}
    I present a geometric multigrid solver for the Poisson equation with continuous cubic elements in three dimensions,
smoothed by a multiplicative vertex-patch Schwarz method with exact local solves. All fields are kept in persistent cell-wise
storage, where every cell holds its own copy of each shared degree of freedom. The patch smoother yields a continuous
solution without any exchange between cells. On the GPU, the patch update is computed in registers, one cell per thread, and the cells of a
patch exchange data by warp shuffles, so shared memory is avoided almost entirely. Despite the low occupancy that this
design implies, the arithmetic of the patch update runs close to the FP32 peak. Fast
diagonalization, merged transforms, and a reflection symmetry of the patch then reduce the arithmetic until a sweep
moves its data at $87\%$ of the bandwidth of a plain device copy. The hierarchy therefore runs in single precision inside a double-precision
outer iteration. Rounding the unassembled residual destroys the cancellation between opposite-sign copies of shared
degrees of freedom; exchanging the copies between cells in double precision performs this cancellation before rounding,
without changing the dual vector they represent. With this \emph{rebalancing}, the cycle converges as in double
precision and independently of the mesh size, the error converges at the optimal rate down to the precision of
the solve, and the solver reaches over $400$
million degrees of freedom per second on a single NVIDIA A100 GPU.

\end{abstract}

\begin{keywords}
    geometric multigrid, vertex-patch smoother, Schwarz methods, fast diagonalization, mixed precision, matrix-free,
    finite elements, GPU
\end{keywords}

\begin{MSCcodes}
    65N55, 65N30, 65F08, 65Y05, 65Y10
\end{MSCcodes}

\section{Introduction}
\label{sec:patch_introduction}

Geometric multigrid with matrix-free, sum-factorized operator evaluation is an established solver for high-order finite
element discretizations of elliptic problems on GPUs~\cite{Kronbichler2012,kronbichler2019multigrid}. The arithmetic
throughput of GPUs has grown faster than their memory bandwidth~\cite{williams2009roofline}, so the operator
evaluations, transfers, and smoothers of a V-cycle run at the rate at which they move their data. Single precision
halves the bytes per value and doubles the arithmetic peak, and since a preconditioner needs no full accuracy, the
hierarchy can run in FP32 inside an FP64 outer iteration without limiting the accuracy of the
solution~\cite{goddeke2007performance,kronbichler2019multigrid,cui2025implementation}.

Halving the traffic shifts the balance toward arithmetic, and a kernel must then also sustain the arithmetic rate.
Matrix-free GPU kernels stage their intermediates in shared memory between the stages of sum
factorization~\cite{ljungkvist2017matrix,kronbichler2019multigrid}, and its bandwidth can become their limit, as for
the GPU patch smoother of~\cite{cui2025implementation}. The register file is both larger and faster than shared
memory~\cite{nvidia2020ampere}, and registers are the direct operands of the arithmetic instructions. The price is occupancy: the more registers a thread
holds, the fewer warps are resident to hide latency by switching between them. A kernel with low occupancy has to hide
latency by independent instructions within each thread instead, and can then still reach the peak arithmetic
rate~\cite{volkov2010better}.

\paragraph{Persistent cell-wise storage}
We build on the storage of~\cite{wichrowski2026DSS,wichrowski2026MG}, on which multigrid runs
without indirect memory access: every field is kept in a redundant element-local representation, and the vector of
unique degrees of freedom is never formed (\secref{sec:patch_setting}). On this storage, a vertex patch assembles its residual from the cell copies
it contains and writes the same correction into every copy, so the smoother preserves continuity by itself: the V-cycle
needs no direct stiffness summation (DSS) or other synchronization, no valence, and no constraint handling, and Dirichlet and refinement-edge degrees
of freedom are never written because no patch reaches them. All global objects are operators applied matrix-free; we speak of operator
application, and matrix language is reserved for the small dense cell and patch matrices. This paper brings mixed
precision, vertex-patch smoothing, and in-register evaluation to multigrid on this storage.

\paragraph{Overlapping Schwarz smoothers}
Schwarz methods on vertex patches, the $2^d$ cells sharing a vertex, are effective smoothers for high-order
discretizations. Pavarino~\cite{pavarino1993additive} proved that additive Schwarz methods with vertex-patch subdomains
are robust in the polynomial degree, Schöberl et al.~\cite{schoberl2008additive} extended this analysis to triangular
and tetrahedral elements, and Arnold, Falk, and Winther~\cite{ArnoldFalkWinther97Hdiv,ArnoldFalkWinther00} established
vertex-patch smoothers for multigrid in $H(\mathrm{div})$ and $H(\mathrm{curl})$. Brubeck and
Farrell~\cite{brubeck2021scalable} made the vertex-star relaxation scalable and robust in the polynomial degree. The
overlapping Schwarz smoothers of spectral element and discontinuous Galerkin methods apply the local inverse by fast
diagonalization~\cite{fischer2000overlapping,LottesFischer05,WitteArndtKanschat21}, and local $p$-multigrid replaces
the exact local solve where the patch problem has no such structure~\cite{wichrowski2025local,wichrowski2026stokespatch}.

\paragraph{Patch smoothers on GPUs}
With exact local solves, a multiplicative vertex-patch smoother gives V-cycles that reduce the residual by orders of
magnitude per cycle~\cite{cui2025implementation}. The price is data movement: in three dimensions every cell belongs to eight
patches, so a sweep moves the data of each cell eight times, where an operator evaluation moves it once
(\secref{sec:patch_smoother}). On a GPU, the
cost of a vertex-patch smoother is therefore set by the bytes it moves and by how much of its arithmetic can be hidden
behind them. Closest to the present work, Cui et al.~\cite{cui2025implementation} realize the multiplicative
vertex-patch smoother on GPUs with each patch held in shared memory, a design bound by shared-memory bandwidth. The
smoother of this paper keeps the patch in registers instead, and its arithmetic runs close to the FP32 peak.

\paragraph{Stagnation and rebalancing}
An FP32 hierarchy on cell-wise storage stalls for a reason absent from assembled vectors. The unassembled residual
carries inter-element fluxes of size $\mathcal{O}(1)$ and of opposite sign, which cancel only on assembly, and rounding
them to FP32 before the assembly loses the small part that carries the residual. With assembled vectors, the cell
contributions are accumulated in FP64 and the fluxes cancel before the residual is rounded. We remove the fluxes by
\emph{rebalancing}, an exchange of the copies between the cells that share a degree of freedom, in FP64, which cancels
the fluxes without changing the dual vector they represent (\secref{sec:patch_mixed}). For CG, one rebalancing lasts
many iterations; the stationary iteration applies it in every cycle.

\paragraph{In-register evaluation}
A thread owns one cell of a patch and holds its $64$ coefficients and all intermediates in registers, the eight cells
of a patch occupy a hypercube of lanes, and all exchange within the patch uses butterfly shuffles; shared memory only
makes the global memory accesses coalesced (\secref{sec:patch_registers}). Applied eight times to data loaded once, a
volumetric $\mathbb{Q}_3$ Laplace kernel built this way keeps the FP32 pipeline $99.6\%$ busy, at $96\%$ of the FP32
peak, despite its low occupancy. A realization that stages through shared memory has no such margin: at $p=3$, the
shared-memory bandwidth is exactly balanced with the FP32 arithmetic it feeds, so any bank conflict, synchronization,
or address computation makes it the bound (\secref{sec:patch_onchip}).

\paragraph{Reducing the arithmetic}
With all intermediates in registers, the time of the smoother is set by its instruction count until it meets the
memory traffic. Accumulator chains, the residual formed in the cell eigenbasis, merged adjacent transforms, and a
reflection symmetry of the patch, whose even and odd modes make the block of the partner cell a signed copy of the
own one, remove instructions without changing the result (\secref{sec:patch_arithmetic}). In a prototype that moves only the
iterate and the correction, the final variant takes $14\%$ longer than a kernel that performs only its memory traffic.

\paragraph{Results}
With rebalancing, the mixed-precision cycle converges as the FP64 one: four V-cycles reduce the residual by nine orders
of magnitude independently of the mesh size, and the error converges at the optimal rate $\mathcal{O}(h^4)$ down to the
precision of the solve (\secref{sec:patch_mixed_convergence}). A complete smoothing sweep moves its data at $87\%$ of the
bandwidth of a plain device copy, and rebalancing takes $5\%$ of the solve time. The solver reaches $416$ million degrees of freedom
per second on one A100, $2.2$ times the throughput reported in~\cite{cui2025implementation}
(\secref{sec:patch_results}).

\section{Problem Setting}
\label{sec:patch_setting}

We consider the Poisson problem on a bounded domain $\Omega\subset\mathbb{R}^3$ with homogeneous Dirichlet boundary
conditions: find $u\in H^1_0(\Omega)$ such that
\begin{equation}
    \int_\Omega \nabla u\cdot\nabla v\,dx=\int_\Omega fv\,dx\qquad\text{for all }v\in H^1_0(\Omega).
    \label{eq:patch_poisson}
\end{equation}
The domain is covered by a mesh $\mathcal{T}$ of hexahedral cells, and~\eqref{eq:patch_poisson} is discretized by
continuous tensor-product elements $\mathbb{Q}_p$ with Gauss--Lobatto nodal bases. The construction below is written
for $p=3$, the degree of all experiments, and for Cartesian meshes.

The restriction to Cartesian meshes is deliberate. Unfitted discretizations such as CutFEM~\cite{burman2015cutfem}
represent the geometry on a Cartesian background mesh, and geometric multigrid with vertex-patch smoothers on that mesh
is effective for them~\cite{cui2025multigrid,wichrowski2026cutfemMG}. Away from the boundary, the patches of the
background mesh are exactly the uniform patches treated here, so their smoother is the kernel of this work. The GPU
patch smoother~\cite{cui2025implementation}, the reference for our comparisons, is likewise restricted to
Cartesian meshes.

On the reference interval let $M,K\in\mathbb{R}^{(p+1)\times(p+1)}$ be the one-dimensional mass and stiffness matrices.
On a cube of edge length $h$ they scale as $sM$ and $s^{-1}K$ with $s=h/2$, so the cell matrix of the Laplacian is
\begin{equation}
    A_K=s\,\bigl(K\otimes M\otimes M+M\otimes K\otimes M+M\otimes M\otimes K\bigr),
    \label{eq:patch_cell_operator}
\end{equation}
and on a uniform mesh level all cells share the same $A_K$. Matrix-free evaluation applies~\eqref{eq:patch_cell_operator}
by sum factorization~\cite{orszag1979spectral,Kronbichler2012}.

All fields are kept in the cell-wise storage of~\cite{wichrowski2026DSS}: every cell holds its own copy of each of its
$(p+1)^3=64$ coefficients, and the vector of unique degrees of freedom (DoFs) is never formed. Let
$\mathbb{V}_{\mathrm{cell}}$ be the space of cell-wise coefficient vectors and $\mathbb{V}_{CG}$ the space of unique
DoFs. The gather $\mathcal{G}\colon\mathbb{V}_{CG}\to\mathbb{V}_{\mathrm{cell}}$ copies each unique coefficient into
all cells that contain it, and its transpose $\mathcal{G}^T$ sums the cell copies. \emph{Primal} vectors (iterates,
corrections) lie in $\operatorname{Im}\mathcal{G}$: their copies agree. \emph{Dual} vectors (right-hand sides,
residuals), marked by a tilde, hold the unassembled cell integrals; their copies carry partial sums that only the
assembly $\mathcal{G}^T$ combines into a functional. The volumetric operator $A$, which applies $A_K$ on
each cell independently, maps primal to dual vectors, and for a primal $z=\mathcal{G}z_{CG}$ and a dual $\tilde r$ the pairing identity
\begin{equation}
    \langle\tilde r,z\rangle=\langle\mathcal{G}^T\tilde r,z_{CG}\rangle
    \label{eq:patch_pairing}
\end{equation}
holds~\cite{wichrowski2026DSS}. A preconditioner maps a dual residual to a primal correction and depends on the
residual only through its assembly, $P=\mathcal{G}P_{CG}\mathcal{G}^T$.

The cell-wise representation is one of two ways to store a continuous field. In the classical one, a vector holds the
unique DoFs; the operator gathers them into the cells, applies the cell matrices, and accumulates the cell results back
into the vector, so every result is assembled~\cite{Kronbichler2012}. In the cell-wise representation the operator
needs no accumulation, but the preconditioner must return a primal vector: one whose copies agree or, more generally,
which satisfies the constraints of the discrete space, such as those of hanging nodes. This is usually enforced by an
exchange of information between cells after each application of the preconditioner, either a gather--scatter through
the vector of unique DoFs, as in Nek5000/NekRS~\cite{nekrs2022,gslib}, or a summation of the copies directly between
neighbouring cells, the direct stiffness summation (DSS)~\cite{deville2002highorder,wichrowski2026DSS}. This work takes
a third way: the preconditioner, a multigrid cycle with the vertex-patch smoother, returns a vector that satisfies the
constraints by construction, and no exchange between cells is needed (\secref{sec:patch_sweep}).

The multigrid hierarchy consists of the nested uniform levels $\mathcal{T}_0\sqsubset\dots\sqsubset\mathcal{T}_L$
obtained by repeated global refinement of a coarse mesh. The V-cycle follows~\cite{wichrowski2026MG}: restriction and
prolongation are the cell-wise transfers of that work, which act on unassembled residuals and on primal corrections
without any DSS, and the residual of each level is evaluated by the matrix-free level operator. The only new component
is the smoother, which replaces the Chebyshev--Jacobi smoother of~\cite{wichrowski2026MG} by the vertex-patch method of
\secref{sec:patch_smoother}.

The outer iteration is either CG preconditioned by the V-cycle or the stationary iteration $u_{k+1}=u_k+M^{-1}\tilde
    r_k$ with $\tilde r_k=\tilde b-Au_k$ and $u_0=0$, where $M^{-1}$ denotes one V-cycle from a zero initial guess. Both
monitor convergence by
\begin{equation}
    \eta_k=\langle\tilde r_k,M^{-1}\tilde r_k\rangle^{1/2}
    \label{eq:patch_monitor}
\end{equation}
relative to $\eta_0$. For the symmetric cycle this is the preconditioned residual norm of CG~\cite{wichrowski2026MG};
for the unsymmetric cycle the pairing is not induced by an inner product and serves as a monitor only. The pairing
requires the V-cycle applied to $\tilde r_k$, so $\eta_k$ is known only after that cycle. In the stationary iteration
this is the cycle that computes the next correction, so the monitor lags one cycle behind the iterate: after $n$
cycles the last known value is $\eta_{n-1}$, and the last correction is never measured. CG applies the V-cycle to its
updated residual before the test, so the last V-cycle of a CG solve serves only the test.

\section{The Vertex-Patch Smoother}
\label{sec:patch_smoother}
\label{sec:patch_sweep}

The smoother is a multiplicative Schwarz method on overlapping vertex patches, the classical choice for high-order
elements~\cite{pavarino1993additive,fischer2000overlapping,LottesFischer05,brubeck2021scalable}. For an interior vertex
$v$ of a level, the vertex patch $\omega_v$ is the set of the $2^3=8$ cells that have $v$ as a corner, and its local
space is spanned by the basis functions of the patch nodes that are not on its boundary. Let $A_v$ be the patch matrix,
the Galerkin restriction of the level operator to this local space, and let $G_{v,K}$ map the patch unknowns to the
coefficients of a cell $K\in\omega_v$; both are given explicitly in \secref{sec:patch_fdm}. Given a primal level
iterate $u$ and the unassembled level right-hand side $\tilde b$, the update of patch $\omega_v$ is
\begin{equation}
    u_K\leftarrow u_K+G_{v,K}\,A_v^{-1}\sum_{K'\in\omega_v}G_{v,K'}^T\bigl(\tilde b_{K'}-A_{K'}u_{K'}\bigr),
    \qquad K\in\omega_v.
    \label{eq:patch_update}
\end{equation}
The sum over $K'$ assembles the patch residual from the raw cell residuals: each patch DoF has copies only in cells of
$\omega_v$, so the sum of $G_{v,K'}^T$ over the patch is the full assembly $\mathcal{G}^T$ for those DoFs. The
update solves the local problem exactly and annihilates the assembled residual on the local space of
$\omega_v$.

One patch update thus consists of four steps, with the residual computed locally on the patch as in~\cite{wichrowski2025smoothers}:
\begin{enumerate}
    \item read $u_K$ and $\tilde b_K$ for the eight cells $K\in\omega_v$;
    \item compute the cell residuals $\tilde b_K-A_Ku_K$ and assemble them into the patch residual;
    \item solve the local problem with $A_v$ and add the correction to every copy of each patch DoF;
    \item write $u_K$ back for the eight cells.
\end{enumerate}

The update writes the same value into every copy of a patch DoF, so a continuous iterate stays continuous. The smoother
therefore needs none of the exchanges between cells described in \secref{sec:patch_setting}: no DSS, no
gather--scatter, no valence of a DoF, and no representation of the continuity constraints. Together with the transfers
and level operators of~\cite{wichrowski2026MG}, which act on unassembled residuals and continuous corrections cell by
cell, this makes the V-cycle free of any DSS. The only exchange between cells in the whole solver is the rebalancing of
\secref{sec:patch_rebalancing}, which serves the precision of the residual and not its continuity.

The patches overlap: an interior DoF of a cell belongs to the patches of all eight corners of that cell. The patches
are therefore partitioned into eight colors by the parities of the integer coordinates of their centre vertices. Two
patches of one color have centres that differ by at least two cells in some direction and share no cell, so the patches
of one color are updated concurrently and without conflicts, and each interior cell belongs to exactly one patch of
each color. The colors are processed one after another, each forming its residual from the iterate left by the previous
one, so a cell is visited eight times per sweep. Reversing the color order gives the adjoint sweep. Forward
pre-smoothing with reverse post-smoothing yields a symmetric V-cycle, which is required when the cycle preconditions
the conjugate gradient method; with forward post-smoothing the cycle is not symmetric and is used as a stationary
iteration.

Only interior vertices carry patches. The Dirichlet DoFs on $\partial\Omega$ lie on the rim of every patch that
contains them and are never written, so the smoother needs no boundary treatment. On a locally refined hierarchy the
same rule, applied to vertices whose eight cells all belong to the level, excludes the DoFs on refinement edges, as
in~\cite{wichrowski2026MG}; the experiments below use uniform hierarchies. The coarsest level of the hierarchy consists
of the eight cells around the centre of the domain, whose single patch is the whole level, so one patch update is an
exact coarse solve.

\subsection{Fast Diagonalization of the Cell-Wise Patch}
\label{sec:patch_fdm}

Along one axis the patch spans two cells and $2p+1=7$ nodes; clamping the two rim nodes leaves $2p-1=5$ nodes per axis
and $125$ unknowns. Let $M_v,K_v\in\mathbb{R}^{5\times5}$ be the one-dimensional mass and stiffness matrices of two
adjacent cells, assembled over the shared node and restricted to the five non-rim nodes. Because all cells of a level
are equal, the patch matrix is the Kronecker sum
\begin{equation}
    A_v=s\,\bigl(K_v\otimes M_v\otimes M_v+M_v\otimes K_v\otimes M_v+M_v\otimes M_v\otimes K_v\bigr)\in\mathbb{R}^{125\times125},
    \label{eq:patch_matrix}
\end{equation}
and its inverse is applied by fast diagonalization~\cite{Lynch1964}. Let $Q\in\mathbb{R}^{5\times5}$ and
$\mu\in\mathbb{R}^5$ solve the generalized eigenproblem $K_vQ=M_vQ\operatorname{diag}(\mu)$ with $Q^TM_vQ=I$. Then
\begin{equation}
    A_v^{-1}=s^{-1}\,(Q\otimes Q\otimes Q)\,\Lambda_v^{-1}\,(Q^T\otimes Q^T\otimes Q^T),
    \qquad (\Lambda_v)_{ijk}=\mu_i+\mu_j+\mu_k .
    \label{eq:patch_fdm}
\end{equation}
Applied as written, this acts on the assembled patch vector of $5^3$ values. In cell-wise storage it would require
three separate steps: summing the copies of the shared nodes across cells into the $125$ unique values, the
transforms of~\eqref{eq:patch_fdm}, and copying the result back into the cells. The five nodes and five modes per axis
also do not divide evenly between the two cells of an axis, so a distribution of the work over the cells would be
unbalanced, with the centre node belonging to neither cell.

In cell-wise storage the patch consists of the copies of its $125$ values in the eight cells. Along one axis, label the
two cells by $c\in\{0,1\}$, the low and the high one. Each cell keeps three of its four nodes, the outer rim node being
clamped: cell~$0$ keeps nodes $1,2,3$ and cell~$1$ nodes $0,1,2$. Let $R_c\in\mathbb{R}^{3 \times4}$ select these
nodes. The two cells together hold $3+3=6$ values per axis, in which the centre node appears twice. The duplication map
\begin{equation}
    g=\begin{bmatrix}g_0\\ g_1\end{bmatrix}=
    \begin{bmatrix}I_3&0\\0&I_3\end{bmatrix}\in\mathbb{R}^{6\times5},
    \qquad g_0,g_1\in\mathbb{R}^{3\times5},
    \label{eq:patch_duplication}
\end{equation}
in which the two identity blocks overlap in the centre column, sends the five patch nodes to these six copies. A cell
$K\in\omega_v$ sits at position $(c_x,c_y,c_z)$ in the patch, and the patch coefficients reach its $64$ coefficients
through
\begin{equation}
    G_{v,K}=R_{c_x}^Tg_{c_x}\otimes R_{c_y}^Tg_{c_y}\otimes R_{c_z}^Tg_{c_z}\in\mathbb{R}^{64\times125}.
    \label{eq:patch_gather}
\end{equation}
The patch thus occupies $27$ of the $64$ coefficients of each of its cells, and $(3+3)^3=216$ cell-wise values in total.

We instead apply fast diagonalization to the cell-wise patch operator $G_vA_v^{-1}G_v^T$, which acts on the $(3+3)^3$
cell-wise values of the patch. Per axis, compose the eigenvectors with the duplication~\eqref{eq:patch_duplication} and
pad by one row,
\begin{equation}
    F=\begin{bmatrix}(gQ)^T\\0\end{bmatrix}\in\mathbb{R}^{6\times6},
    \qquad \tilde\mu=(\mu_1,\dots,\mu_5,1),
    \label{eq:patch_fdm_cellwise_1d}
\end{equation}
and let $\tilde\Lambda_v$ be the $6^3$ diagonal with entries $\tilde\mu_i+\tilde\mu_j+\tilde\mu_k$. Since the padded
row of $F$ is zero, every modal coefficient that involves the padded mode vanishes, its eigenvalue only has to be
nonzero, and
\begin{equation}
    G_vA_v^{-1}G_v^T=s^{-1}\,(F^T\otimes F^T\otimes F^T)\,\tilde\Lambda_v^{-1}\,(F\otimes F\otimes F),
    \qquad G_v=g\otimes g\otimes g .
    \label{eq:patch_fdm_cellwise}
\end{equation}
The forward transform $F$ maps the six cell-wise values of an axis to modal coefficients and sums the two copies of the
centre node in the same contraction; the backward transform $F^T$ produces the correction directly in both copies.
The assembly and the scatter are part of the transforms, and the $125$ unique values are never formed.

The padding makes the modal dimension equal to the cell-wise one, six, so that both split evenly between the two cells.
Assign three modes to each cell and partition $F$ into $3\times3$ blocks $F_{ab}$, where $a$ is the cell that owns the
modes and $b$ the cell that owns the nodes. Along one axis the forward and backward transforms of cell $a$, with
partner $\bar a=1-a$, are
\begin{equation}
    y_a=F_{aa}x_a+F_{a\bar a}x_{\bar a},
    \qquad
    x_a=F_{aa}^Ty_a+F_{\bar aa}^Ty_{\bar a}.
    \label{eq:patch_fdm_blocks}
\end{equation}
In three dimensions each of the eight cells holds a $3\times3\times3$ share of the data at every stage: its nodal values
before the forward transform and an octant of the $6^3$ modal coefficients after it. Each cell performs the same
operations on the same amount of data, which is what the GPU realization of \secref{sec:patch_registers} maps onto
eight threads. Compared with~\eqref{eq:patch_fdm}, the transforms have six instead of five rows; the reductions of
\secref{sec:patch_arithmetic} more than recover this cost.

\section{Mixed Precision}
\label{sec:patch_mixed}

The multigrid hierarchy is a preconditioner, and its components need not be more accurate than the reduction it
achieves per cycle. Running it in FP32 has two benefits. Every vector of the hierarchy takes half the memory and half
the traffic, which matters most for the smoother: one sweep visits every cell eight times, each visit reading $u$ and
$\tilde b$ and writing $u$. And on current GPUs the FP32 vector units have twice the peak of the FP64 ones, so
arithmetic that is close to the memory traffic in FP64 moves well below it in FP32. Both gains are realized only if the
cycle converges as fast as in FP64.

All computations in this work run on one NVIDIA A100-SXM4-80GB GPU~\cite{nvidia2020ampere}, whose peak rates are
listed in Table~\ref{tab:patch_a100}. TF32 is a tensor-core format with the exponent range of FP32 and $11$
significand bits.

\begin{table}[htbp]
    \centering
    \begin{tabular}{@{}l r@{}}
        \toprule
        Memory bandwidth                          & $2039$\,GB/s       \\
        Shared-memory bandwidth                   & $19.5$\,TB/s       \\
        FP64, vector units                        & $9.7$\,TFLOP/s     \\
        FP32, vector units                        & $19.5$\,TFLOP/s    \\
        FP64, tensor cores                        & $19.5$\,TFLOP/s    \\
        TF32, tensor cores                        & $156$\,TFLOP/s     \\
        \bottomrule
    \end{tabular}
    \caption{Peak rates of the NVIDIA A100-SXM4-80GB~\cite{nvidia2020ampere}.}
    \label{tab:patch_a100}
\end{table}

\subsection{Precision Split}
\label{sec:patch_precision_split}

The outer iteration, either the conjugate gradient method or a stationary iteration of V-cycles, runs in FP64: the
iterate, the right-hand side, the operator on the finest mesh, the residual, and all inner products. The V-cycle is an
FP32 inner solve, in line with established mixed-precision multigrid on
GPUs~\cite{goddeke2007performance,kronbichler2019multigrid,cui2025implementation}. On entry the FP64 residual is
rounded to FP32 and copied into the level hierarchy; on exit the FP32 correction is converted back and added to the
FP64 iterate. The correction holds identical values in all copies of a DoF, its conversion to FP64 is exact, and the
accumulation in FP64 adds identical values to identical copies, so the FP64 iterate stays continuous. All level vectors, the patch smoother, the level operators, and the transfers are FP32. The patch
smoother uses plain FP32 fused multiply-adds.

The level operators, which form the residuals restricted between the levels, and the transfers are sum-factorization
kernels whose one-dimensional contractions run on the tensor cores in TF32. To recover FP32 accuracy, both operands are
split into a TF32-exact high part and a remainder, $x=x_h+x_l$ and $A=A_h+A_l$, and the product is evaluated as
\begin{equation}
    Ax\approx A_hx_h+\bigl(A_lx_h+A_hx_l\bigr),
    \label{eq:patch_tf32_split}
\end{equation}
dropping the term $A_lx_l$ of relative size $2^{-22}$, with all terms accumulated in FP32. For $p=3$ the contracted
dimension is four, while the smallest TF32 tile has a contraction dimension of eight, so the stacked operand $[x_h\mid
            x_l]$ fills a tile that a single unsplit operand would leave half empty. The split operators reach a relative error of
about $10^{-6}$ against FP64, the accuracy of FP32 arithmetic in the same sum factorization.

\subsection{Test Configuration}
\label{sec:patch_test_configuration}

All experiments solve~\eqref{eq:patch_poisson} on the unit cube with the manufactured solution
$u(x)=\prod_{d=1}^3\sin(\pi x_d)$ and $f=3\pi^2u$, discretized by $\mathbb{Q}_3$ elements. The finest mesh is obtained
by $L$ global refinements of the unit cube and has $8^L$ cells and $(3\cdot2^L+1)^3$ unique DoFs, including the
Dirichlet boundary. DoF counts and throughputs refer to these unique DoFs. The hierarchy runs from the level of eight cells to the finest level. Each
V-cycle applies one pre- and one post-smoothing sweep.

Kernel times are CUDA-event averages over repeated launches after warm-up. Solve times are wall-clock times of the
complete solve, synchronized with the device, excluding the setup of the hierarchy and one warm-up cycle.

The tolerance is $\eta/\eta_0\le10^{-9}$ for the monitor~\eqref{eq:patch_monitor}, and the $L^2$ error is reported as
an independent check. For the stationary iteration we extrapolate the unmeasured last correction at the rate of the
last two cycles,
\begin{equation}
    \hat\eta_n=\eta_{n-1}^2/\eta_{n-2},
    \label{eq:patch_extrapolation}
\end{equation}
and report the average contraction $\rho=(\hat\eta_n/\eta_0)^{1/n}$ and the
fractional iteration count
\begin{equation}
    n_{\mathrm{frac}}=\frac{\log 10^{-9}}{\log\rho},
    \label{eq:patch_fractional_iterations}
\end{equation}
the number of cycles at the average contraction needed to reach the tolerance. For CG, which measures after each
update, $\rho$ and $n_{\mathrm{frac}}$ use the final measured residual. The fractional count distinguishes solvers whose
integer counts coincide.

\subsection{Stagnation}
\label{sec:patch_cancellation}

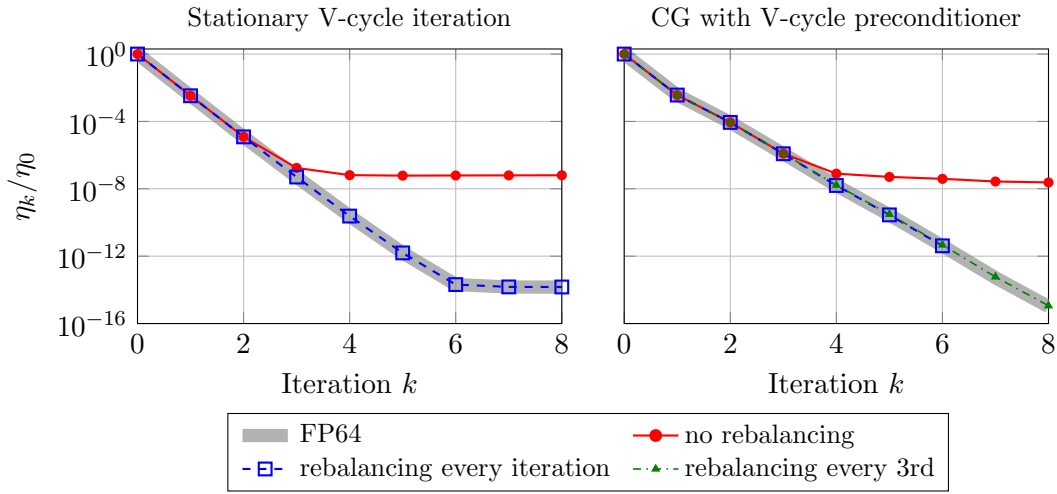
\begin{figure}[htbp]
    \centering
    \pgfplotstableread[col sep=comma]{raw_results/patch_mixed/precision_history_L6/precision_history_L6.csv}\precisionhistory
\pgfplotsset{
    patch history axis/.style={
            scale only axis,
            width=0.34\linewidth,
            height=0.22\linewidth,
            xlabel={Iteration $k$},
            xmin=0, xmax=8,
            xtick={0,2,...,8},
            ymin=1e-16, ymax=2,
            ytick={1,1e-4,1e-8,1e-12,1e-16},
            grid=major,
            title style={font=\small, yshift=-3pt},
            mark size=1.5pt,
            unbounded coords=jump,
            filter discard warning=false,
            restrict y to domain=-36.9:1,
        },
    patch history fp64/.style={gray, line width=5pt, opacity=0.6},
    patch history none/.style={red, thick, mark=*},
    patch history every/.style={blue, thick, dashed, mark=square, mark size=2.4pt, mark options={solid}},
    patch history third/.style={green!50!black, semithick, dashdotted, mark=triangle*, mark size=1.4pt, mark options={solid}},
}
\begin{tikzpicture}
    \begin{semilogyaxis}[
            patch history axis,
            name=cycleaxis,
            title={Stationary V-cycle iteration},
            ylabel={$\eta_k/\eta_0$},
        ]
        \addplot[patch history fp64] table[x=k, y=fp64_cycle] {\precisionhistory};
        \addplot[patch history none] table[x=k, y=cycle_rebalance0] {\precisionhistory};
        \addplot[patch history every] table[x=k, y=cycle_rebalance1] {\precisionhistory};
    \end{semilogyaxis}
    \begin{semilogyaxis}[
            patch history axis,
            name=cgaxis,
            at={(cycleaxis.south east)}, anchor=south west, xshift=0.05\linewidth,
            title={CG with V-cycle preconditioner},
            yticklabels={},
            legend to name=patchhistorylegend,
            legend columns=2,
            legend cell align=left,
            legend style={font=\small, /tikz/every even column/.append style={column sep=6pt}},
        ]
        \addplot[patch history fp64] table[x=k, y=fp64_cg_15it] {\precisionhistory};
        \addlegendentry{FP64}
        \addplot[patch history none] table[x=k, y=cg_rebalance0] {\precisionhistory};
        \addlegendentry{no rebalancing}
        \addplot[patch history every] table[x=k, y=cg_rebalance1] {\precisionhistory};
        \addlegendentry{rebalancing every iteration}
        \addplot[patch history third] table[x=k, y=cg_rebalance3_15it] {\precisionhistory};
        \addlegendentry{rebalancing every 3rd}
    \end{semilogyaxis}
    \path ($(cgaxis.east) + (cycleaxis.west) - (current bounding box.west)$);
\end{tikzpicture}

\smallskip
\pgfplotslegendfromname{patchhistorylegend}
    \caption{Residual history of the mixed-precision iteration on the mesh with $L=6$, without and with rebalancing
        (\secref{sec:patch_rebalancing}), against the same iteration in FP64 (grey bands). Left: stationary iteration
        with the forward V-cycle; right: CG preconditioned by the symmetric V-cycle. CG with rebalancing in every iteration terminates at $k=6$ on
        an absolute threshold of the residual.}
    \label{fig:patch_precision_history}
\end{figure}

Figure~\ref{fig:patch_precision_history} compares the residual history of the mixed-precision iteration on the mesh
with $L=6$ with that of the same iteration in FP64. The first two iterations agree with and without the rebalancing of
\secref{sec:patch_rebalancing}; from the third on, the iteration without it falls behind the FP64 one and stalls. The
stationary cycle stops at $\eta/\eta_0\approx6\cdot10^{-8}$ from the fourth cycle on, and CG fluctuates near
$2\cdot10^{-8}$. In a fixed number of five CG iterations the effect shows in the discretization error: it stagnates
near $10^{-8}$ from $L=6$ on (Table~\ref{tab:patch_norebalance}).

The cause is the rounding of the unassembled residual. The copy of the residual that cell $K$ holds for a DoF with
basis function $\varphi_i$ is, after integration by parts,
\begin{equation}
    \tilde r_{K,i}=\int_K f\varphi_i\,dx-\int_K\nabla u_h\cdot\nabla\varphi_i\,dx
    =\int_K(f+\Delta u_h)\varphi_i\,dx-\int_{\partial K}\partial_{n_K}u_h\,\varphi_i\,ds .
    \label{eq:patch_residual_flux}
\end{equation}
For a DoF on a face shared by $K$ and $K'$, the boundary terms of the two copies are fluxes of the discrete solution
through the common face with opposite normals. Their size is set by the gradient of the solution, $\mathcal{O}(1)$
relative to the data, and it does not decrease as the iteration converges. The assembly $\mathcal{G}^T$ adds the copies,
the two fluxes cancel up to the jump of the normal derivative, and what remains is the residual that drives the
iteration. Once the iteration has converged somewhat, the copies therefore have the form $a+e$ and $-a+e'$ with
$|e|,|e'|\ll|a|$, and the assembled residual is $e+e'$.

Summed in FP64, the copies give $e+e'$ to FP64 accuracy. If they are rounded to FP32 first, each carries an error of up
to $\varepsilon_{32}|a|$; the subsequent assembly in the patch smoother still cancels $a$, but the rounding errors of
the two copies do not cancel, and they are no longer small compared with $e+e'$. The $\mathcal{O}(1)$ flux survives the
rounding and is then removed by the assembly, while the small part that carries the residual is lost. Let
$\varepsilon_{32}=2^{-24}$ be the unit roundoff and $\mathrm{fl}$ the rounding to FP32. Since $\mathcal{G}^T$ sums at
most $n_{\max}=8$ copies, $\|\mathcal{G}^T\|_2=\sqrt{n_{\max}}$ and
\begin{equation}
    \frac{\|\mathcal{G}^T(\mathrm{fl}(\tilde r)-\tilde r)\|}{\|\mathcal{G}^T\tilde r\|}
    \le\sqrt{n_{\max}}\;\varepsilon_{32}\;\kappa(\tilde r),
    \qquad
    \kappa(\tilde r)=\frac{\|\tilde r\|}{\|\mathcal{G}^T\tilde r\|}.
    \label{eq:patch_cancellation_bound}
\end{equation}
The cancellation ratio $\kappa$ grows as the residual decreases while the fluxes keep their size. When it
reaches $\varepsilon_{32}^{-1}$, the assembled residual that the V-cycle sees is rounding noise, and the outer iteration
stalls. A random vector shows no
cancellation, so the effect is visible only on residuals of an actual solve.

A related issue arises at the Dirichlet boundary. The raw values of $Au$ on the Dirichlet DoFs are of the size of the
operator times the solution, much larger than the residual, and are discarded by the solver. They are set to zero in
the FP64 residual before it is rounded; otherwise they set the scale of the FP32 cycle and stall it at a relative
residual of about $10^{-8}$.

\subsection{Rebalancing}
\label{sec:patch_rebalancing}

The remedy is to cancel the fluxes before rounding, by exchanging the copies of each DoF between the cells that share
it, in FP64. This replaces $a+e$ and $-a+e'$ by $(e+e')/2$ each: the flux is gone, and the rounding error becomes relative to the residual itself. Let
$D=\mathcal{G}^T\mathcal{G}$, which multiplies each unique DoF by its valence, and
\begin{equation}
    \Pi=\mathcal{G}D^{-1}\mathcal{G}^T\colon\mathbb{V}_{\mathrm{cell}}^*\to\mathbb{V}_{\mathrm{cell}}^*,
    \label{eq:patch_rebalance}
\end{equation}
which replaces every copy by the mean of all copies of its DoF: in practice, the exchange is an averaging. $\Pi$ is a
symmetric projection, and $\Pi\tilde r$ represents the same dual vector as $\tilde r$,
\begin{equation}
    \mathcal{G}^T\Pi\tilde r=\mathcal{G}^T\mathcal{G}D^{-1}\mathcal{G}^T\tilde r=\mathcal{G}^T\tilde r .
    \label{eq:patch_rebalance_assembly}
\end{equation}
The patch smoother, which sums the copies itself, therefore needs no modification. Its cancellation ratio is bounded independently of the
state of the iteration: with $\|D^{-1}\|_2\le1$,
\begin{equation}
    \kappa(\Pi\tilde r)=\frac{\|\mathcal{G}D^{-1}\mathcal{G}^T\tilde r\|}{\|\mathcal{G}^T\tilde r\|}\le\sqrt{n_{\max}},
    \label{eq:patch_rebalance_kappa}
\end{equation}
so by~\eqref{eq:patch_cancellation_bound} the rounding perturbs the assembled residual by at most
$n_{\max}\varepsilon_{32}$ relative to its size, however far the iteration has progressed.

$\Pi$ is computed by the dimension-split DSS of~\cite{wichrowski2026DSS}. The DSS processes the three axes in turn. In
the pass along one axis, every DoF on a face normal to that axis has one copy in each of the two cells that share the
face, and the pass replaces both copies $a$ and $b$ by $a+b$; DoFs on edges and vertices, shared by four and eight
cells, are summed completely after two and three passes. For $\Pi$ each pass writes $(a+b)/2$ instead. The factor
$1/2$ leaves the represented vector unchanged: the two copies still sum to $a+b$, so~\eqref{eq:patch_rebalance_assembly} holds
after every pass, whereas the plain sum would double it. The passes compose into the mean over all
$2^k$ copies of a DoF shared by $2^k$ cells, so no valence needs to be stored, and one rebalancing costs one FP64 DSS
on the finest mesh.

For any primal $z=\mathcal{G}z_{CG}$, the pairing identity~\eqref{eq:patch_pairing}
and~\eqref{eq:patch_rebalance_assembly} give $\langle\Pi\tilde r,z\rangle=\langle\tilde r,z\rangle$. In exact
arithmetic, replacing the residual of CG by $\Pi\tilde r$ in place therefore changes neither its inner products nor the
preconditioned residual $P\tilde r=\mathcal{G}P_{CG}\mathcal{G}^T\tilde r$, so the CG iterates are unchanged. CG then
updates its residual recursively, $\tilde r\leftarrow\tilde r-\alpha A p$, and each update reintroduces cancellation
only of the size of the current step $\alpha Ap$, which shrinks as CG converges. One rebalancing therefore lasts many
iterations, and it suffices to rebalance every few iterations. A stationary iteration, in contrast, forms its residual
anew as $\tilde b-Au$ in every cycle. That residual carries the full opposite-sign fluxes of $\tilde b$ and $Au$, and
it is rebalanced in every cycle.

\subsection{Convergence with Rebalancing}
\label{sec:patch_mixed_convergence}

With rebalancing, the mixed-precision iteration follows the FP64 one past the former stall
(Figure~\ref{fig:patch_precision_history}); the two histories agree to two or three digits in every iteration. The
stationary cycle keeps its contraction of about $3\cdot10^{-3}$ per cycle until it levels off at
$\eta/\eta_0\approx1.5\cdot10^{-14}$ from the seventh cycle on, the level at which the FP64 cycle levels off as well.
CG with rebalancing in every third iteration follows FP64 CG down to $10^{-15}$ at iteration $8$, and rebalancing in
every iteration gives the same history.

Table~\ref{tab:patch_convergence} shows the stationary iteration with forward post-smoothing and rebalancing in every
cycle. Four cycles reach the tolerance on every mesh, and the fractional count is virtually identical on all meshes
and to that of the same cycle in FP64: rebalancing recovers the convergence of the FP64 cycle. The $L^2$ error
converges at the optimal rate up to $L=6$. From $L=7$ on, it reaches the precision of the solve: the extrapolated
discretization error at $L=8$, about $4\cdot10^{-12}$, is below the algebraic error that the tolerance of $10^{-9}$
leaves, and the rate drops. At $L=7$, CG, which overshoots the tolerance further, keeps the rate close to the
optimal~$4$.

\begin{table}[htbp]
    \centering
    \begin{tabular}{@{}r r r r r r r r r@{}}
    \toprule
    & & & \multicolumn{3}{c}{mixed FP64/FP32} & FP64 & & \\
    \cmidrule(lr){4-6}\cmidrule(lr){7-7}
    $L$ & cells & DoFs & $n$ & $n_{\mathrm{frac}}$ & $\rho$ & $n_{\mathrm{frac}}$ & $\|u-u_h\|_{L^2}$ & rate \\
    \midrule
    4 & $16^3$  & $1.18\cdot10^{5}$ & 4 & 3.70 & 0.0037 & 3.70 & $2.419\cdot10^{-7}$  & --   \\
    5 & $32^3$  & $9.13\cdot10^{5}$ & 4 & 3.73 & 0.0038 & 3.73 & $1.514\cdot10^{-8}$  & 4.00 \\
    6 & $64^3$  & $7.19\cdot10^{6}$ & 4 & 3.72 & 0.0038 & 3.72 & $9.470\cdot10^{-10}$ & 4.00 \\
    7 & $128^3$ & $5.71\cdot10^{7}$ & 4 & 3.71 & 0.0037 & 3.71 & $6.879\cdot10^{-11}$ & 3.78 \\
    8 & $256^3$ & $4.55\cdot10^{8}$ & 4 & 3.70 & 0.0037 & ---  & $3.244\cdot10^{-11}$ & 1.08 \\
    \bottomrule
\end{tabular}

    \caption{Stationary V-cycle iteration, forward post-smoothing, tolerance $10^{-9}$. Mixed precision with
        rebalancing in every cycle against the same cycle with all components in FP64. $n$ is the number of cycles,
        $n_{\mathrm{frac}}$ the fractional count~\eqref{eq:patch_fractional_iterations}, $\rho$ the average
        contraction. The FP64 cycle was not run for $L=8$.}
    \label{tab:patch_convergence}
\end{table}

Table~\ref{tab:patch_norebalance} shows what the lost cancellation costs at a fixed iteration count. Both runs apply
five CG iterations with the mixed-precision cycle. Without
rebalancing the $L^2$ error stagnates near $10^{-8}$: from $L=6$ on, the rounding of the residual dominates the
discretization error and the convergence rate is lost. With rebalancing every third iteration, the error converges at
the optimal rate.

\begin{table}[htbp]
    \centering
    \begin{tabular}{@{}r r r r r@{}}
    \toprule
    & \multicolumn{2}{c}{rebalancing every 3rd iteration} & \multicolumn{2}{c}{no rebalancing} \\
    \cmidrule(lr){2-3}\cmidrule(lr){4-5}
    $L$ & $\|u-u_h\|_{L^2}$ & rate & $\|u-u_h\|_{L^2}$ & rate \\
    \midrule
    4 & $2.419\cdot10^{-7}$  & --   & $2.421\cdot10^{-7}$ & --   \\
    5 & $1.514\cdot10^{-8}$  & 4.00 & $1.624\cdot10^{-8}$ & 3.90 \\
    6 & $9.468\cdot10^{-10}$ & 4.00 & $1.157\cdot10^{-8}$ & 0.49 \\
    7 & $6.221\cdot10^{-11}$ & 3.93 & $6.697\cdot10^{-9}$ & 0.79 \\
    \bottomrule
\end{tabular}

    \caption{Five CG iterations with the mixed-precision V-cycle, with rebalancing in every third iteration
        and without rebalancing.}
    \label{tab:patch_norebalance}
\end{table}

\begin{remark}
    Mixed-precision multigrid is often combined with GMRES~\cite{cui2025implementation}. On cell-wise storage, GMRES
    needs an inner product on the primal Krylov space, and its residual estimate rests on an Arnoldi relation that an
    FP32 V-cycle satisfies only to single precision. Since a Krylov method gains less than one cycle over the stationary
    iteration (\secref{sec:patch_results_outer}), we leave a mixed-precision GMRES for cell-wise storage to future
    work.
\end{remark}

\section{In-Register Evaluation}
\label{sec:patch_registers}

A GPU kernel keeps its working data on chip in two places. The register file is the larger and the faster one: each
streaming multiprocessor of the A100 has $256$\,KB of registers against at most $164$\,KB of shared
memory~\cite{nvidia2020ampere},\footnote{The H100 (Hopper) has the same $256$\,KB of registers and up to $228$\,KB of
shared memory per streaming multiprocessor~\cite{nvidia2022hopper}, and so does the Blackwell B200~\cite{nvidia2025blackwell}.
The B200 additionally has $256$\,KB of \emph{tensor memory} per streaming multiprocessor, which holds the accumulator
and optionally the first operand of a tensor-core operation, the second being read from shared memory, and thus lowers the register pressure of tensor-core kernels~\cite{nvidia2025ptx}.} and registers are the direct operands of the arithmetic instructions. A fused multiply-add reads three of them, so at the FP32
peak the register operands amount to about $117$\,TB/s, six times the shared-memory bandwidth. Shared memory, in turn, is addressable by all threads of a block, but every access is a
separate load or store instruction through the memory pipeline, accesses of a warp that fall into the same bank are
serialized, and data written by one thread is visible to another only after a barrier.

Registers are allocated by the compiler, not by the program, in two steps. The kernel is first translated into an
intermediate assembly (PTX) that assumes an unlimited number of virtual registers, one per value; only the translation
to machine code (SASS) assigns them to physical registers, at most $255$ per thread or fewer if the programmer sets a
bound for the kernel. Registers are private to a thread and are addressed only by indices known at compile time, so an
array kept in registers requires all loops over it to be unrolled; an array indexed at run time is placed in local
memory. Values that do not fit into the physical registers are \emph{spilled} to local memory, a per-thread memory
backed by DRAM and cached in L1 and L2, whose traffic competes with that of the kernel's data. The register count also
sets the \emph{occupancy}: the $65\,536$ registers of a streaming multiprocessor are divided among its resident
threads, so a kernel with $255$ registers per thread admits $8$ of the $64$ possible warps. Fewer warps hide less
latency by switching between them; a kernel with low occupancy has to hide it by independent instructions within each
thread~\cite{volkov2010better}. Limiting the register count to raise the occupancy trades latency for spill traffic.
For a matrix-free operator limited by cache traffic, Böhm et al.~\cite{bohm2026performance} report that a register
bound admitting more warps caused spills and reduced the throughput to $58\%$ and $30\%$ of the spill-free kernel.

In practice, matrix-free GPU kernels keep their intermediates in shared memory. Sum factorization contracts a
coefficient cube along each axis in turn, and consecutive contractions need the data distributed differently across
threads; since a thread cannot address the registers of another, the cube is staged through shared memory between
stages~\cite{ljungkvist2017matrix,kronbichler2019multigrid}. Shared memory consequently becomes the resource around
which the kernels are designed: the number of elements per thread block is chosen to fit a shared-memory
budget~\cite{prusak2026mixed}, DoFs shared between cells are staged in shared memory for a tile of
cells~\cite{bohm2026performance}, and the GPU vertex-patch smoother~\cite{cui2025implementation} holds
the whole patch there. The only exchange of registers between threads is the warp shuffle, which passes one value per
lane and instruction among the $32$ lanes of a warp.

The kernel of this work evaluates the complete patch update in registers. A thread owns one cell of a patch and holds
all its coefficients, every intermediate of the residual and of the patch solve stays in the registers of the thread
that computes it, and the coupling between the cells of a patch moves by warp shuffles. Each field is moved between
DRAM and the registers once per cell visit; shared memory is used only to make these global accesses coalesced
(\secref{sec:patch_layout}). The price is a register count close to the limit and the corresponding low occupancy.

To test whether the arithmetic can reach the peak under these conditions, we apply the volumetric $\mathbb{Q}_3$
Laplace kernel eight times to data loaded once, which effectively removes the memory bandwidth from the measurement
(Table~\ref{tab:patch_laplace_peak}). The kernel follows the same design: one cell per thread, the $64$ coefficients
and every intermediate of the sum factorization in registers, and no shared memory. On $2^{19}$ cells, the busiest
pipeline of the streaming multiprocessors runs at $98.7\%$ of its throughput. The remaining gap comes from the size of
the launch: $2^{19}$ cells are too few to keep all streaming multiprocessors busy until its end. On $2^{22}$ cells the
pipeline utilization rises to $99.6\%$. Hence, the low occupancy does not prevent the kernel from saturating the FP32
units.

\begin{table}[htbp]
    \centering
    \begin{tabular}{@{}r r r r r r r@{}}
    \toprule
    cells & time [ms] & TFLOP/s & FP32 peak & SM pipeline & DRAM & L1 \\
    \midrule
    $2^{19}$ & 0.839 & 18.56 & $95.2\%$ & $98.7\%$ & $14.8\%$ & $4.0\%$ \\
    $2^{22}$ & 6.647 & 18.74 & $96.1\%$ & $99.6\%$ & $15.7\%$ & $4.1\%$ \\
    \bottomrule
\end{tabular}

    \caption{Volumetric $\mathbb{Q}_3$ Laplace kernel applied eight times per load, FP32, A100, one cell per thread
        with all intermediates in registers. Each application costs $3712$ FLOP per cell, $58$ FLOP per byte of DRAM
        traffic over the eight applications. Time and rates from Nsight Compute; \emph{SM pipeline} is the utilization
        of the busiest pipeline of the streaming multiprocessors, \emph{DRAM} the fraction of the memory bandwidth, and
        \emph{L1} the throughput utilization of the unit that serves the L1 cache, shared memory, and register
        spills. Its load comes from the global memory accesses alone: neither run spills registers or accesses shared
        memory.}
    \label{tab:patch_laplace_peak}
\end{table}

\subsection{Thread Mapping}
\label{sec:patch_thread_mapping}

The unit of work of a thread is a cell. A contraction along any axis of that cell is a sequence of fused multiply-adds
on the thread's own registers, and no stage requires a redistribution of data between threads. Only the coupling
between the cells of a patch crosses threads, and the cell-wise fast diagonalization of \secref{sec:patch_fdm} reduces
it to one partner per axis. This partition is natural for cell-wise storage, where a cell is the unit in which data is
stored, and it is balanced because the $(3+3)^3$ formulation gives every cell the same share of the patch. It differs
from the conventional partition by DoFs of the assembled patch, in which every contraction stage couples values held by
different threads and the data is redistributed through shared memory after each stage (\secref{sec:patch_onchip}).
Table~\ref{tab:patch_mapping} summarizes the two designs.

\begin{table}[htbp]
    \centering
    \begin{tabular}{@{}l l l@{}}
        \toprule
                                 & DoF-wise partition~\cite{cui2025implementation}               & Cell-wise partition (this work)    \\
        \midrule
        Patch representation     & assembled, $7^3$ values                                     & cell-wise, $8\times4^3$ values     \\
        Unit of work of a thread & one DoF or one column of DoFs                               & one cell                           \\
        Threads per patch        & $49$ to $343$                                               & $8$                                \\
        Intermediates held in    & shared memory                                               & registers                          \\
        Exchange between stages  & shared memory, block barrier                                & warp shuffle with one partner      \\
        Shared memory used for   & all stages                                                  & coalescing global accesses only    \\
        Limiting resource        & shared-memory bandwidth                                     & DRAM bandwidth or FP32 units$^{*}$ \\
        \bottomrule
    \end{tabular}
    \caption{Two ways of distributing the update of one $\mathbb{Q}_3$ vertex patch over GPU threads.
        $^{*}$The limit depends on the work per loaded cell. For the constant-coefficient Laplacian of this work, one
        update per load is bound by DRAM bandwidth (Table~\ref{tab:patch_stages}); with more arithmetic per load, as with
        several updates per load, it is bound by the FP32 units, which it drives at up to $86\%$ of their peak (\secref{sec:patch_arithmetic_rate}).}
    \label{tab:patch_mapping}
\end{table}

One color of the sweep is one kernel launch, and each warp processes four patches. The eight cells of a patch occupy
eight consecutive lanes, ordered as
\begin{equation}
    \text{lane}=8q+4c_z+2c_y+c_x,
    \label{eq:patch_lanes}
\end{equation}
for the $q$-th patch of the warp. The three bits of the lane index within a patch are the position of the cell along
the three axes, so the eight lanes form a three-dimensional hypercube with the same geometry as the patch. The partner
of a cell across axis $d$ is the lane whose index differs in bit $d$, and the exchange with it is a single butterfly
shuffle with lane mask $2^d$, which exchanges register values between the two lanes without shared memory and without
synchronization beyond the implicit one of the shuffle. Each contraction stage of~\eqref{eq:patch_fdm_blocks} along axis
$d$ pairs every lane with its partner across $d$; the three axes use the three dimensions of the hypercube in turn,
as in a radix-two fast Fourier transform.

The block structure~\eqref{eq:patch_fdm_blocks} assigns each lane the same amount of work at every stage. In the
forward stages, a lane receives the partner's values of each pencil by shuffle and produces the three modes it owns
along the contracted axis; after the stage along $x$ it holds three $x$-modes for the $4\times4$ $(y,z)$ pencils of its
cell, after $y$ a $3\times3\times4$ block, and after $z$ a $3\times3\times3$ octant of the $6^3$ modal cube. The modal
cube is thus distributed over the lanes exactly as the cells of the patch are, and the diagonal scaling by
$\tilde\Lambda_v^{-1}$ needs no communication: each lane selects the eigenvalues of its octant from its lane bits. The
backward stages retrace the exchange pattern in reverse order, $z$, $y$, $x$, and each lane ends with the $27$
correction values of its own cell.

Per patch update and thread, the forward stages shuffle the four values of each pencil, $64$, $48$, and $36$ values
along $x$, $y$, and $z$ as the cube shrinks, and the three backward stages shuffle the three values of each of the $9$
pencils of the $3\times3\times3$ block, $81$ values. These $229$ shuffles, against $1776$ floating-point
instructions in the prototype of Table~\ref{tab:patch_stages}, are the entire communication within a patch.

All lanes execute the same instruction stream. What differs between them is data, selected once at the start of the
kernel from the three lane bits: the $3\times4$ forward blocks and $3\times3$ backward blocks of the modes the lane
owns along each axis, the signs that the reflection symmetry of \secref{sec:patch_parity} attaches to the partner's
values, the eigenvalues of the lane's modal octant, and the position of the clamped rim node of the cell along each
axis. Each lane loads only its own blocks, about $90$ values in total, and keeps them in registers. There is no
divergence between the lanes of a patch.

A shuffle requires all lanes of the warp to participate. If the number of patches of a color is not a multiple of four,
the lanes of the missing patches therefore compute on a copy of the first patch and skip only their stores.

The complete color step holds the coefficients of $u$ and of $\tilde b$ during the transforms and uses about $250$
registers per thread without spilling, which admits eight resident warps per streaming multiprocessor. Latency is
hidden by the independent arithmetic on the $64$ coefficients of each thread rather than by many resident warps.
Capping the register count to admit more warps causes spilling and is slower: $128$ and $96$ registers per thread
increase the time of a color step by factors of $1.9$ and $3.2$.

\subsection{Memory Layout}
\label{sec:patch_layout}

The solver uses two layouts. The FP64 vectors of the outer iteration on the finest mesh, on which the residual is
evaluated and rebalanced, are stored in the blocked layout of~\cite{wichrowski2026DSS}, with $32$ cells per block and
the cell index as the fastest axis; this is the layout in which the DSS that realizes the rebalancing is
efficient~\cite{wichrowski2026MG}. The FP32 level vectors of the V-cycle need neither DSS nor rebalancing and are
stored cell-contiguously, with the $64$ coefficients of a cell adjacent, so that a whole cell is fetched at once and a
patch gathers its eight cells by index. The transfer of the residual into the hierarchy narrows it to FP32 and changes
the layout in one pass, and the correction is widened and converted back in the same way.

With one cell per thread in the cell-contiguous layout, a warp's loads of one coefficient are $256$ bytes apart, each
touching a separate sector, and a direct realization in this layout runs $2.9$ times slower than in a blocked layout in
which the cell index is the fastest axis (Table~\ref{tab:patch_layout}). The loss comes from the access pattern, not
from the layout. A vector access, of type \texttt{float4} in CUDA, moves four adjacent FP32 words, $16$ aligned bytes,
between memory and four registers of one thread in a single instruction. It reduces the number of instructions, but
coalescing depends only on the sectors touched by the warp: with one cell per thread, vector loads touch the same
scattered sectors and gain $6\%$. The kernel therefore separates the access pattern from the thread mapping. Each thread
reads the index of its own cell from the patch table, and obtains the indices of the cells it helps to load by shuffle
from their owners. In each round, eight consecutive lanes load the eight $16$-byte pieces of one $128$-byte row of one
cell by vector loads, so that one load instruction of the warp reads four whole rows, and the values are handed to
their owning threads through shared memory. The staging buffer holds half a cell per thread
and round: staging whole cells needs $8.7$\,KB per warp, which makes shared memory instead of registers the occupancy
limit and costs $11\%$. All loads of a field are issued before the first hand-over, so the warp waits on memory once
per field, and strides of $36$ and $27$ words place concurrent accesses on distinct banks.

The correction occupies $27$ of the $64$ coefficients of a cell, at positions that depend on the corner of the cell in
its patch. Storing them directly would be a strided store. Instead the correction is expanded to the full cube with
zeros at the clamped rim positions, staged through shared memory, and added to a second read of the cell of $u$, so
that every store writes whole rows. The second read is served by the L2 cache. Its loads are all issued before the
first store, so their latencies overlap instead of adding up, which improves the color step by $8\%$ to $23\%$. The DRAM traffic per cell visit is thus $3\cdot256$ bytes for
$u$, $\tilde b$, and the store, plus $4$ bytes of the cell index.

\begin{table}[htbp]
    \centering
    \begin{tabular}{@{}l r r r@{}}
    \toprule
    Layout and access & kernel [ms] & stub [ms] & stub [GB/s] \\
    \midrule
    Blocked, $32$ cells per block                 & 0.257 & 0.226 & 1686 \\
    Cell-contiguous, one cell per thread          & 0.752 & 0.795 & 480  \\
    Cell-contiguous, staged through shared memory & 0.268 & 0.219 & 1742 \\
    \bottomrule
\end{tabular}

    \caption{One color step with blocked and with cell-contiguous storage, $2^{20}$ cells, FP32, one application per
        load. The kernel is variant~5 of the prototype of Table~\ref{tab:patch_stages}. The stub performs the loads and stores of the kernel without arithmetic.}
    \label{tab:patch_layout}
\end{table}

In the final kernel the loads are hidden behind the arithmetic: removing the second read of $u$ and the store together
reduces a color step on $2^{23}$ cells from $4.06$ to $3.06$\,ms, and keeping all cell data in cache instead of DRAM
reduces it only further to $2.95$\,ms. The remaining difference to the arithmetic-only time is the store.

\subsection{Relation to Shared-Memory Realizations}
\label{sec:patch_onchip}

The GPU vertex-patch smoother of Cui et al.~\cite{cui2025implementation} is the shared-memory design applied to the
patch: the patch is held in its assembled numbering, a thread is assigned to one DoF or, for higher degrees in three
dimensions, to one column of $2p+1$ DoFs, one or several patches form a thread block, and every stage of the residual
evaluation and of the fast diagonalization reads its operands from shared memory and writes its results back, separated
by block-wide synchronizations. After removing all bank conflicts, their $\mathbb{Q}_3$ kernel in three dimensions runs
at $95\%$ of the throughput of the on-chip memory and at $24\%$ of the FP64 peak. Their on-chip roofline confirms
shared-memory bandwidth as the limit for all degrees. The highest fraction, $36\%$ of the peak, is reached at the
higher degrees $5$ to $7$, and equally in FP64 and in FP32~\cite{cui2025implementation}: halving the word size does not
lift the kernel off its bound.

At $p=3$ this bound is not a matter of tuning: the shared-memory traffic of a staged contraction is as large as the
arithmetic it serves. A contraction stage of a $4\times4\times4$ cube with $4\times4$ matrices performs $4$ fused
multiply-adds per entry and, staged through shared memory, writes and reads each entry once: $8$ bytes for $4$ fused
multiply-adds in FP32. At the FP32 peak of $9.75\cdot10^{12}$ fused multiply-adds per second, the shared-memory bandwidth of
Table~\ref{tab:patch_a100} delivers likewise $2$ bytes per fused multiply-add. The shared-memory bandwidth and the FP32 peak are
thus exactly balanced for this kernel, and any bank conflict, synchronization, or address computation makes shared
memory the bound. Moreover, the transformations of \secref{sec:patch_arithmetic} reduce the arithmetic per entry to two
terms per mode while the traffic per stage stays the same; with shared-memory staging, they would make the kernel
slower relative to its bound rather than faster.

The same balance appears in tensor-core realizations of the cell operator, where the axis rotations between
contractions are lowered through shared memory. In our implementation of the $\mathbb{Q}_3$ cell operator on the A100
tensor cores, these rotations take $33\%$ of the run time in FP64 and $44\%$ to $65\%$ in the reduced-precision
variants. However well shared memory is used, a kernel that stages its intermediates there competes for a resource
whose bandwidth matches the arithmetic it feeds.

The in-register kernel avoids this bound instead of approaching it. Its shared-memory traffic is fixed by the data it
moves, not by its arithmetic: each of $u$, $\tilde b$, and the correction passes through shared memory once, written by
one thread and read by another, $1536$ bytes per cell visit or about twice the DRAM traffic. At the DRAM rate of a
color step on $2^{23}$ cells (Table~\ref{tab:patch_sweep}) this is $3.2$\,TB/s, about $16\%$ of the $19.5$\,TB/s
shared-memory bandwidth of the device. Table~\ref{tab:patch_shmem} compares the measured usage with the shared-memory
realization~\cite{cui2025implementation}.

\begin{table}[htbp]
    \centering
    \begin{tabular}{@{}l r r r@{}}
    \toprule
    & \cite{cui2025implementation} & \multicolumn{2}{c}{this work} \\
    \cmidrule(lr){2-2}\cmidrule(lr){3-4}
    & conflict-free kernel & prototype & final kernel \\
    \midrule
    Precision                                   & FP64   & FP32     & FP32    \\
    Shared memory used in the arithmetic        & yes    & no       & no      \\
    Shared-memory bytes per cell visit          & --     & 0        & 1536    \\
    L1/shared-memory throughput (of peak)       & 95.3\%  & 42\%     & 45\%    \\
    Shared-memory bank conflicts                & 0      & 0        & 4.2\%   \\
    MIO throttle stall                          & --     & 2.0\%    & 6.3\%   \\
    Arithmetic throughput (of peak)             & 24\%$^{*}$ & 60\%  & 41\%    \\
    \bottomrule
\end{tabular}

    \caption{Shared-memory usage of vertex-patch kernels for $\mathbb{Q}_3$ in three dimensions on the A100. The data
        of~\cite{cui2025implementation} are those of its conflict-free fused kernel. The prototype is variant~5 of
        Table~\ref{tab:patch_stages}; the final kernel is the cell-contiguous color step of
        Table~\ref{tab:patch_sweep}. Bytes per cell visit of the final kernel are counted from its staging. The MIO
        throttle stall is the fraction of warp cycles stalled because the instruction queue of the memory input/output
        (MIO) unit, which serves shared-memory accesses and warp shuffles, is full.
        $^{*}$At degrees $5$ to $7$ the kernel of~\cite{cui2025implementation} reaches $36\%$ of the peak, in FP64 and in FP32 alike.}
    \label{tab:patch_shmem}
\end{table}

The measurement matches the count. The final kernel moves $12.5$ shared-memory wavefronts of $128$ bytes per cell
visit, $1600$ bytes, which is the counted $1536$ bytes plus $4.2\%$ of bank conflicts.\footnote{We did not pursue
these conflicts, since the kernel is close to its memory traffic anyway (Table~\ref{tab:patch_stages}).} The L1 throughput includes the
global loads, which pass through the same unit: the prototype, which uses no shared memory, reaches $42\%$ from them
alone. Neither kernel is limited by the on-chip memory, and the final kernel, bound by DRAM bandwidth, leaves $59\%$ of
the FP32 peak unused.

\section{Reducing the Arithmetic}
\label{sec:patch_arithmetic}

With all intermediates in registers, the time of a color step is set by its instruction count until it meets the memory
traffic. This section removes instructions in four steps that do not change the result. Table~\ref{tab:patch_stages}
records their effect, measured in a prototype that stores cells in the blocked layout of~\cite{wichrowski2026DSS}. The
prototype reads the $64$ coefficients of $u$ and writes the $27$ values of the correction.
It omits the right-hand side, whose transform costs the final kernel three more contractions, and the
read-modify-write of \secref{sec:patch_layout}. The reference row is a traffic stub that performs the same loads and stores without arithmetic and bounds what any kernel
with this data movement can achieve.

\begin{table}[htbp]
    \centering
    \begin{tabular}{@{}l r r r r r@{}}
    \toprule
    Organization of the arithmetic & time [ms] & GB/s & \% of stub & FP instr. & rel.\ error \\
    \midrule
    1. Nodal contractions & 0.526 & 725  & 43 & 3615 & $2.4\cdot10^{-7}$ \\
    2. One accumulator chain & 0.468 & 816  & 48 & 3129 & $2.1\cdot10^{-7}$ \\
    3. Residual in the cell eigenbasis           & 0.430 & 887  & 53 & 2873 & $3.0\cdot10^{-7}$ \\
    4. Merged adjacent transforms                & 0.379 & 1007 & 60 & 2423 & $2.5\cdot10^{-7}$ \\
    5. Reflection symmetry                       & 0.259 & 1472 & 88 & 1776 & $2.2\cdot10^{-7}$ \\
    \midrule
    Traffic stub, no arithmetic                  & 0.227 & 1680 & 100 & -- & -- \\
    \bottomrule
\end{tabular}

    \caption{One color step of the vertex-patch smoother on $2^{20}$ cells ($131\,072$ patches), FP32, one application
        per load, A100. Each row computes the same result as the previous one. Instruction counts are per thread in the
        compiled code. The bandwidth counts the field data the step must move. The relative error is measured against
        an FP64 reference that applies the Kronecker sum and the dense patch inverse directly.}
    \label{tab:patch_stages}
\end{table}

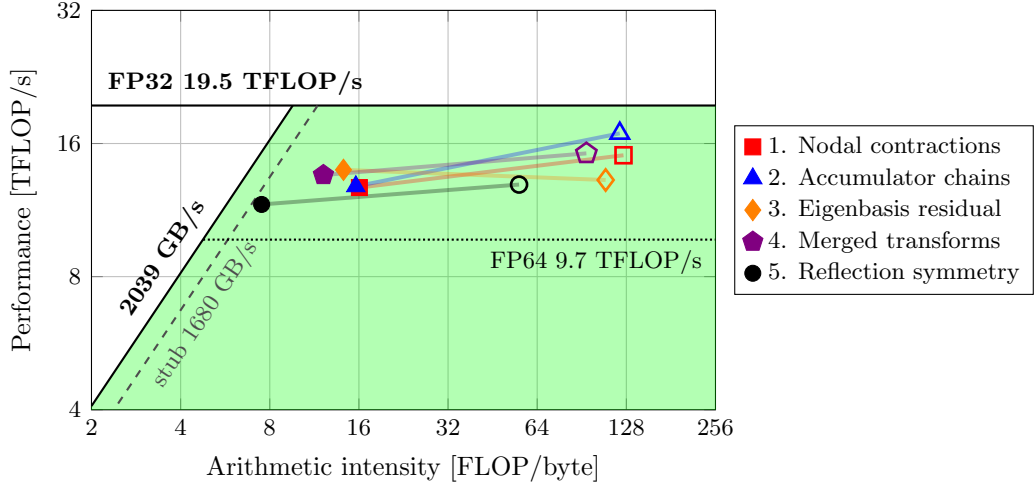
\begin{figure}[htbp]
    \centering
    \begin{tikzpicture}
    \begin{loglogaxis}[
            scale only axis,
            width=0.5\linewidth,
            height=0.32\linewidth,
            xlabel={Arithmetic intensity [FLOP/byte]},
            ylabel={Performance [TFLOP/s]},
            log basis x=2,
            xmin=2, xmax=256,
            xtick={2,4,8,16,32,64,128,256},
            xticklabels={2,4,8,16,32,64,128,256},
            log basis y=2,
            ymin=4, ymax=32,
            ytick={4,8,16,32},
            yticklabels={4,8,16,32},
            grid=both,
            legend style={at={(1.03,0.5)}, anchor=west},
            legend columns=1,
            legend cell align=left,
            legend style={font=\footnotesize},
            every axis label/.append style={font=\small},
            tick label style={font=\footnotesize},
        ]
        \addplot[green, fill=green, draw=none, fill opacity=0.3, forget plot]
        coordinates {(0.001,0.002039) (9.56,19.5) (1000,19.5) (1000,0.0001) (0.001,0.0001)};

        \addplot[thick, domain=2:9.56, samples=2, forget plot] {2.039*x}
        node[pos=0.45, sloped, above, font=\footnotesize\bfseries] {2039 GB/s};
        \addplot[thick, domain=2:256, samples=2, forget plot] {19.5};
        \addplot[thick, dashed, black!70, domain=2:11.6, samples=2, forget plot] {1.680*x}
        node[pos=0.42, sloped, below, font=\footnotesize] {stub 1680 GB/s};
        \addplot[thick, densely dotted, domain=4.76:256, samples=2, forget plot] {9.7};

        \node[above right, font=\footnotesize\bfseries] at (axis cs:2.1,19.5) {FP32 19.5 TFLOP/s};
        \node[below left, font=\footnotesize] at (axis cs:250,9.7) {FP64 9.7 TFLOP/s};

        \addplot[red, line width=1.5pt, opacity=0.3, forget plot] coordinates {(16.06,12.72) (125.2,15.05)};
        \addplot[blue, line width=1.5pt, opacity=0.3, forget plot] coordinates {(15.62,12.8) (121.63,16.87)};
        \addplot[orange, line width=1.5pt, opacity=0.3, forget plot] coordinates {(14.19,13.95) (108.98,13.23)};
        \addplot[violet, line width=1.5pt, opacity=0.3, forget plot] coordinates {(12.14,13.58) (93.81,15.2)};
        \addplot[black, line width=1.5pt, opacity=0.3, forget plot] coordinates {(7.51,11.65) (55.59,12.92)};

        \addplot[red, mark=square*, mark size=3pt, only marks] coordinates {(16.06,12.72)};
        \addlegendentry{1. Nodal contractions}
        \addplot[blue, mark=triangle*, mark size=4pt, only marks] coordinates {(15.62,12.8)};
        \addlegendentry{2. Accumulator chains}
        \addplot[orange, mark=diamond*, mark size=4pt, only marks] coordinates {(14.19,13.95)};
        \addlegendentry{3. Eigenbasis residual}
        \addplot[violet, mark=pentagon*, mark size=4pt, only marks] coordinates {(12.14,13.58)};
        \addlegendentry{4. Merged transforms}
        \addplot[black, mark=*, mark size=3pt, only marks] coordinates {(7.51,11.65)};
        \addlegendentry{5. Reflection symmetry}

        \addplot[red, mark=square, mark size=3pt, mark options={line width=1pt}, only marks, forget plot]
        coordinates {(125.2,15.05)};
        \addplot[blue, mark=triangle, mark size=4pt, mark options={line width=1pt}, only marks, forget plot]
        coordinates {(121.63,16.87)};
        \addplot[orange, mark=diamond, mark size=4pt, mark options={line width=1pt}, only marks, forget plot]
        coordinates {(108.98,13.23)};
        \addplot[violet, mark=pentagon, mark size=4pt, mark options={line width=1pt}, only marks, forget plot]
        coordinates {(93.81,15.2)};
        \addplot[black, mark=o, mark size=3pt, mark options={line width=1pt}, only marks, forget plot]
        coordinates {(55.59,12.92)};
    \end{loglogaxis}
\end{tikzpicture}
    \caption{Roofline of the color step variants of Table~\ref{tab:patch_stages}, FP32, A100, $2^{20}$ cells. Filled
        markers: one update per load of the cell data, the regime of the smoother. Open markers: eight updates per
        load, which removes the field traffic and shows the arithmetic rate; for variants $3$ and $4$ they include
        register spills. Operations are counted from the executed FP32 instructions, a fused multiply-add as two, over
        the profiled kernel time. The dashed line is the bandwidth of a stub with the same loads and stores and no
        arithmetic, and the dotted line the FP64 peak of the vector units, for reference.}
    \label{fig:patch_roofline}
\end{figure}

\subsection{Accumulator Chains}
\label{sec:patch_accumulator}

The first change does not alter the arithmetic. Written as a sum of products of own and partner terms, each inner loop
was compiled into separate multiply, fused multiply-add, and add instructions, nine instructions for twelve
floating-point operations. Accumulating every product into a single chain lets every term after the first become one
fused multiply-add, which removes $488$ instructions per thread and $11\%$ of the run time.

\subsection{Residual in the Cell Eigenbasis}
\label{sec:patch_eigen_residual}

The cell operator admits the same factorization as the patch. Let $S\in\mathbb{R}^{4\times4}$ and
$\lambda\in\mathbb{R}^4$ solve $KS=MS\operatorname{diag}(\lambda)$ with $S^TMS=I$, so that $M=S^{-T}S^{-1}$ and
$K=S^{-T}\operatorname{diag}(\lambda)S^{-1}$. Then $A_K=s\,(S^{-T})^{\otimes3}\Lambda_K(S^{-1})^{\otimes3}$ with
$(\Lambda_K)_{ijk}=\lambda_i+\lambda_j+\lambda_k$, and the cell residual is
\begin{equation}
    \tilde b_K-A_Ku_K=(S^{-T})^{\otimes3}\,\hat r_K,\qquad
    \hat r_K=(S^{T})^{\otimes3}\tilde b_K-s\,\Lambda_K\,(S^{-1})^{\otimes3}u_K .
    \label{eq:patch_residual_eigen}
\end{equation}
Six contractions and a pointwise product replace the seven contractions of the direct tensor-product evaluation, and
each contraction rotates the coefficient cube by one axis, so three of them restore the original axis order without
the four explicit transposes of the nodal form. For cubic elements $\operatorname{cond}(S)=2.94$, so the change of
basis costs no accuracy in single precision.

\subsection{Merged Transforms}
\label{sec:patch_merged}

In the patch update the backward transform $(S^{-T})^{\otimes3}$ of~\eqref{eq:patch_residual_eigen} is followed
immediately by the selection of the kept nodes and the forward transform of~\eqref{eq:patch_fdm_cellwise}. The three
are multiplied once in setup into
\begin{equation}
    B_{ab}=F_{ab}R_bS^{-T}\in\mathbb{R}^{3\times4},\qquad a,b\in\{0,1\},
    \label{eq:patch_merged}
\end{equation}
so the forward stage along one axis reads $y_a=B_{aa}\hat x_a+B_{a\bar a}\hat x_{\bar a}$ on the eigen coefficients
$\hat x$ of the cell residual. The nodal residual is never materialized, the selection of kept nodes needs no
instructions, and the coefficient cube shrinks along the contracted axis,
$4\times4\times4\to4\times4\times3\to4\times3\times3\to3\times3\times3$. Collecting the factors, one patch update is
\begin{equation}
    \begin{split}
        \delta_K & =s^{-1}\bigl(R_{c_x}^TF_{:c_x}^T\otimes R_{c_y}^TF_{:c_y}^T\otimes R_{c_z}^TF_{:c_z}^T\bigr)\,
        \tilde\Lambda_v^{-1}\sum_{K'\in\omega_v}\bigl(B_{:c'_x}\otimes B_{:c'_y}\otimes B_{:c'_z}\bigr)\,\hat r_{K'}, \\
        u_K      & \leftarrow u_K+\delta_K,
    \end{split}
    \label{eq:patch_update_factored}
\end{equation}
where $F_{:b}\in\mathbb{R}^{6\times3}$ and $B_{:b}\in\mathbb{R}^{6\times4}$ denote the column blocks of cell $b$. This
merge removes three contractions and is the largest single gain in Table~\ref{tab:patch_stages}. It exists only because
two fast diagonalizations are applied back to back and nothing requires the nodal residual between them.

\subsection{Reflection Symmetry}
\label{sec:patch_parity}

The Gauss--Lobatto nodes are symmetric, and the two cells along each axis of a patch have the same size, as on every
level of a uniformly refined hierarchy. The reflection of a patch axis about the centre vertex, which maps patch
node $j$ to $4-j$, therefore commutes with $M_v$ and $K_v$, and the reflection $n\mapsto 3-n$ of a cell axis commutes with $M$
and $K$. The eigenvectors therefore have definite parity: $Q$ has three even and two odd columns, and the columns of
$S$ alternate between even and odd. The reflection also exchanges the two cells of the axis and maps kept node $t$ of
one cell to kept node $2-t$ of the other. We order the modes by parity and assign the three even modes to cell~$0$, and
the two odd modes with the padded one to cell~$1$, so that all modes owned by cell $a$ have one parity $\sigma_a=\pm1$.
Let $\tau\in\{\pm1\}^4$ hold the parities of the cell modes and $J_3$ be the $3\times3$ reversal. Then
\begin{equation}
    B_{a\bar a}=\sigma_a\,B_{aa}\operatorname{diag}(\tau),
    \qquad
    F_{\bar aa}^T=\sigma_{\bar a}\,J_3F_{\bar a\bar a}^T.
    \label{eq:patch_parity}
\end{equation}

The partner's block is a signed copy of the own one. In the forward stages,
\begin{equation}
    y_a=B_{aa}\bigl(\hat x_a+\sigma_a\,\tau\circ\hat x_{\bar a}\bigr),
    \label{eq:patch_fold}
\end{equation}
where $\circ$ is the entrywise product: the partner's four shuffled coefficients are folded into the own ones with one
signed addition each, and one $3\times4$ block is contracted instead of two. In the backward stages,
\begin{equation}
    x_a=F_{aa}^Ty_a+\sigma_{\bar a}\,J_3\bigl(F_{\bar a\bar a}^Ty_{\bar a}\bigr),
    \label{eq:patch_mirror}
\end{equation}
and the bracket is the result the partner computes for its own nodes anyway. Each lane contracts only its own block,
shuffles the result, and adds the partner's result in mirrored order with one signed fused multiply-add per value; the
off-diagonal blocks are never contracted. In the cell operator, the alternating parity of $S$ means that each mode needs
only the sum or only the difference of the two mirror pairs of a pencil: a butterfly of four operations, then two terms
per mode instead of four. The instruction count of every stage was predicted before it was
implemented; the compiled final stage has exactly the predicted $1773$ floating-point instructions per thread, plus
three selects.

\subsection{Arithmetic Rate}
\label{sec:patch_arithmetic_rate}

Figure~\ref{fig:patch_roofline} places the variants on the roofline of the device. Because no operand leaves the
registers, the arithmetic runs close to the peak of the FP32 units. With the data held fixed and the update applied
eight times per load, so that only arithmetic remains, the accumulator-chain variant executes $16.9$\,TFLOP/s, $86\%$
of the FP32 peak, at $92\%$ of the speed of light of the streaming multiprocessors, and the nodal
variant $15.1$\,TFLOP/s. The stall that counts the shuffle and shared-memory pipeline stays at or below $2\%$ in all
variants. Even at one update per load, the four variants without the reflection symmetry lie right of the ridge point:
their arithmetic, not the memory traffic, sets their time, and they reach $65\%$ to $72\%$ of the FP32 peak.

At eight updates per load, variants $2$ to $4$ reach the limit of $255$ registers per thread and spill. Variant $2$
spills about $0.2$\,KB per cell, which the L1 cache absorbs, so its rate is that of the arithmetic. Variants $3$ and
$4$ spill $3.9$ and $2.0$\,KB per cell, and part of the spill stores of variant $3$ reaches DRAM; their open markers
therefore include spill traffic, and variant $3$ executes fewer operations per second at eight updates per load than at
one. The spills are not fatal: the loss of variant $3$ is small, and all spilling variants stay well above the FP64
peak of the vector units, close to the rates of the spill-free ones (Figure~\ref{fig:patch_roofline}).

All variants, at one update per load as well as at eight, run above the $9.7$\,TFLOP/s peak of the FP64 vector units.
The same kernel in FP64 would therefore be slower than its FP32 counterpart even with unlimited memory bandwidth: the
arithmetic rate itself, not only the halved traffic, requires FP32. The FP64 tensor cores of the A100 match the FP32
peak, but the arithmetic of the patch update does not map onto them well. Their instruction multiplies an $8\times4$ by
a $4\times8$ tile distributed over the warp in a fixed fragment layout. The contractions here are $3\times4$ blocks
applied to pencils of a cell held by one thread, so using the tile requires padding three rows to eight, moving the
cell data from the thread that owns it into the fragment layout and back for every stage, and a common matrix for all
lanes, whereas the blocks of the patch update depend on the position of the cell. For the cell operator alone, without
the patch solve, our tensor-core realization reaches $11.6$\,TFLOP/s in FP64 at eight applications per load, with a
third of its time spent on rearranging data between contractions; the register-resident FP32 variants exceed this rate
for the complete patch update.

The reflection symmetry needs cells of equal size within a patch; without it, variants $1$ to $4$ still apply and run
at two thirds of the FP32 peak. Variable coefficients break fast diagonalization itself, and local
$p$-multigrid~\cite{wichrowski2025local} then replaces the exact local solve.

For the constant-coefficient Laplacian, the reductions shift the final variant left of the ridge point. It executes
fewer operations per second, $12.9$\,TFLOP/s at eight updates per load, and is $1.7$ times faster: the transformations
remove work rather than speed up its execution. Time per update, not the fraction of peak, is therefore the measure we
use. At one update per load the final variant is within $14\%$ of the traffic stub; the DRAM, at $74\%$ of its speed of
light, and the streaming multiprocessors, at $85\%$, are then both close to saturation.

\section{Performance}
\label{sec:patch_results}

We now assemble the components into the complete solver: the final smoother kernel in an FP32 V-cycle with
rebalancing, inside an FP64 outer iteration. We measure the throughput of the smoother sweep, choose the outer
iteration, and time the solve. All runs use the configuration of \secref{sec:patch_test_configuration}, and the
convergence of the mixed-precision cycle is that of \secref{sec:patch_mixed_convergence}.

\subsection{Smoother Throughput}
\label{sec:patch_results_smoother}

Table~\ref{tab:patch_sweep} reports the final kernel. A random order of the patches within a color, the worst case
for the gather, is slower; the setup sorts each color by storage order, which avoids it. A complete sweep of eight
colors moves its data at $1521$\,GB/s, $87\%$ of the $1756$\,GB/s that a plain device-to-device copy reaches on
the same GPU. The correction agrees with an FP64 evaluation of the same color step to FP32 accuracy.

\begin{table}[htbp]
    \centering
    \begin{tabular}{@{}l r r r r@{}}
    \toprule
    Case & cells & cell visits & time [ms] & GB/s \\
    \midrule
    One color, patches in storage order & $2^{20}$ & $2^{20}$ & 0.540 & 1499 \\
    One color, random patch order       & $2^{20}$ & $2^{20}$ & 0.591 & 1370 \\
    One color, patches in storage order & $2^{23}$ & $2^{23}$ & 4.085 & 1585 \\
    One color, random patch order       & $2^{23}$ & $2^{23}$ & 4.495 & 1441 \\
    \midrule
    Full sweep, eight colors            & $128^3$  & $16\,387\,064$ & 8.316 & 1521 \\
    \bottomrule
\end{tabular}

    \caption{Final cell-contiguous kernel, FP32. The bandwidth counts $772$ bytes per cell visit: $u$ and $\tilde b$
        read, $u$ written, and the cell index. The second read of $u$ in the read-modify-write is served by the L2
        cache and not counted. The full sweep visits $8\cdot127^3$ cells, eight per interior patch.}
    \label{tab:patch_sweep}
\end{table}

\subsection{Choice of the Outer Iteration}
\label{sec:patch_results_outer}

Table~\ref{tab:patch_solvers} compares outer iterations with the same smoother. The unsymmetric cycle as a stationary
iteration is the fastest: $3.7$ cycles, against $4.7$ CG iterations with the symmetric cycle, each of which also costs
one V-cycle, and $7.2$ cycles for the symmetric cycle without acceleration. As a stationary iteration the symmetric
cycle has a contraction factor fifteen times larger than that of the unsymmetric one. A Krylov method
cannot recover much from the unsymmetric cycle either: the best combination of the corrections it returns, which bounds
what any Krylov method over those vectors can achieve, reduces the residual after a given number of cycles only by a
factor of two to five, less than one cycle.

\begin{table}[htbp]
    \centering
    \begin{tabular}{@{}l l r r r r r@{}}
    \toprule
    Outer iteration & post-smoother & rebalancing & $L=4$ & $L=5$ & $L=6$ & $L=7$ \\
    \midrule
    V-cycle & forward & every cycle       & 3.70 & 3.73 & 3.72 & 3.71 \\
    V-cycle & reverse & every cycle       & 7.25 & 7.24 & 7.24 & 7.23 \\
    CG      & reverse & every iteration   & 4.77 & 4.73 & 4.71 & 4.65 \\
    CG      & reverse & every 3rd iteration & 4.78 & 4.74 & 4.73 & 4.68 \\
    \bottomrule
\end{tabular}

    \caption{Fractional iteration counts~\eqref{eq:patch_fractional_iterations} of outer iterations with the same
        mixed-precision V-cycle, tolerance $10^{-9}$, on meshes with $L$ global refinements. The integer counts are $4$
        for the forward cycle, $8$ for the reverse cycle, and $5$ for CG on every mesh.}
    \label{tab:patch_solvers}
\end{table}

The forward cycle, one pre- and post-smoothing sweep with exact local solves and an unsymmetric cycle,
matches the configuration of~\cite{cui2025implementation}, which reports three iterations for $\mathbb{Q}_3$ in three
dimensions against our four cycles. The two counts describe the same amount of work, as the following remark explains.

\begin{remark}[Counting V-cycles]
    \label{rem:patch_cycle_count}
    An iteration of GMRES is not one application of the preconditioner: $m$ steps apply it $m+1$ times, once to the
    initial residual and once per Arnoldi step, with left preconditioning and, unless the preconditioned basis vectors
    are stored, with right preconditioning as well. The timings in~\cite[Table~7]{cui2025implementation} are for GMRES preconditioned by a V-cycle, so their three iterations at $\mathbb{Q}_3$ are four V-cycles,
    the same number as in Table~\ref{tab:patch_convergence}. The two stopping criteria differ: theirs is the Euclidean
    norm of the assembled residual relative to $\|b\|$, ours $\eta/\eta_0$.
\end{remark}

\subsection{Solve Time}
\label{sec:patch_results_time}

Table~\ref{tab:patch_solve_time} reports the solve time of the stationary iteration of
Table~\ref{tab:patch_convergence}. The timed solves apply a fixed number of four cycles, the count that
Table~\ref{tab:patch_convergence} establishes for the tolerance on every mesh. A solve that stops on the measured
monitor needs a fifth cycle, since the monitor lags one cycle behind the iterate (\secref{sec:patch_setting}). Without
a previous run, the extrapolation~\eqref{eq:patch_extrapolation} estimates the monitor of the current iterate after
every cycle, and stopping on $\hat\eta_n$ avoids the extra cycle; $\hat\eta_n$ is an estimate, not a bound. On the finest mesh the solver reaches $416$ million DoFs per
second.

Table~\ref{tab:patch_profile} breaks down the GPU time of the solve on the mesh with $L=6$. The patch smoother takes
$60\%$ of it, so the reduction of its arithmetic in \secref{sec:patch_arithmetic} acts on the dominant part of the
solve. Rebalancing costs $5\%$, an acceptable price for recovering the convergence of the FP64 cycle.

\begin{table}[htbp]
    \centering
    \begin{tabular}{@{}l r r@{}}
    \toprule
    Part & time [ms] & share \\
    \midrule
    Patch smoother, all levels                      & 12.93 & 60.2\% \\
    Transfers, including FP64/FP32 conversion       & 2.71  & 12.6\% \\
    Elementwise vector operations                   & 1.79  & 8.3\%  \\
    Cell operators, FP64 and FP32                   & 1.22  & 5.7\%  \\
    Rebalancing (DSS)                               & 1.12  & 5.2\%  \\
    Zero fills                                      & 0.96  & 4.5\%  \\
    Convergence monitor $\langle\tilde r,M^{-1}\tilde r\rangle$ & 0.70  & 3.3\%  \\
    Other, device idle between launches             & 0.05  & 0.2\%  \\
    \midrule
    Total (wall time)                               & 21.48 & 100\% \\
    \bottomrule
\end{tabular}

    \caption{GPU time of one solve of Table~\ref{tab:patch_solve_time} on the mesh with $L=6$, four cycles, average of
        ten solves. Shares are relative to the wall time; \emph{Other} is the time in which no kernel runs. The cell operators of the FP64 residual and of the FP32 levels share one kernel and are not
        separated.}
    \label{tab:patch_profile}
\end{table}

On small meshes the solve time is set by kernel launches rather than by work. Every level of a V-cycle issues sixteen
launches for its two smoothing sweeps, one per color, in addition to its level operator and transfers, and this number
does not depend on the mesh size. The solves for $L=4$ and $L=5$ take nearly the same time, although the second mesh has eight
times as many cells. In addition, the colors of a small level do not contain enough patches to
saturate the device. The throughput levels off on the two finest meshes.

\begin{table}[htbp]
    \centering
    \begin{tabular}{@{}r r r r r r@{}}
    \toprule
    $L$ & cells & DoFs & cycles & solve [ms] & MDoF/s \\
    \midrule
    4 & $16^3$  & $1.18\cdot10^{5}$ & 4 & 8.6   & 13.7  \\
    5 & $32^3$  & $9.13\cdot10^{5}$ & 4 & 8.7   & 104.9 \\
    6 & $64^3$  & $7.19\cdot10^{6}$ & 4 & 22.1  & 324.8 \\
    7 & $128^3$ & $5.71\cdot10^{7}$ & 4 & 142.2 & 401.4 \\
    8 & $256^3$ & $4.55\cdot10^{8}$ & 4 & 1093.2 & 416.0 \\
    \bottomrule
\end{tabular}

    \caption{Solve time of the stationary forward V-cycle iteration with rebalancing in every cycle, tolerance
        $10^{-9}$, A100.}
    \label{tab:patch_solve_time}
\end{table}

For comparison, Cui et al.~\cite{cui2025implementation} solve the same problem with $\mathbb{Q}_3$ elements on the same GPU model
with GMRES in FP64 preconditioned by an FP32 V-cycle, at three iterations, i.e.\ four V-cycles on the finest level
(Remark~\ref{rem:patch_cycle_count}), at $188$ million DoFs per second on the mesh with $L=8$, counting unique DoFs as
we do. On this mesh, our solver reaches $2.2$ times their throughput. Both solves apply four V-cycles, which carry most of the work;
they differ only in the vector operations of GMRES and in the stopping criterion.

\section{Conclusions}
\label{sec:patch_conclusion}
We presented a geometric multigrid solver for continuous $\mathbb{Q}_3$ elements on persistent cell-wise storage,
smoothed by a multiplicative vertex-patch Schwarz method with exact local solves. Each patch assembles its residual from
the cell copies and writes the same correction into every copy, so the V-cycle keeps the iterate continuous without any
exchange between cells, valence, or constraint handling. Diagonalizing the cell-wise rather than the assembled patch
absorbs assembly and scatter into the transforms and splits the local solve into eight identical shares, one per cell.

An FP32 hierarchy stalls for a reason unique to the cell-wise formulation: the unassembled residual carries
inter-element fluxes of opposite sign, and rounding them before the assembly loses the small part that carries the
residual. Exchanging the copies between the cells that share a DoF, in FP64 and before rounding, cancels the fluxes
without changing the dual vector they represent. With this rebalancing, the mixed-precision cycle converges as the FP64 one,
four cycles reduce the residual by nine orders of magnitude on every mesh, and the error converges at the optimal rate down to the
precision of the solve.
This is the only exchange between cells in the solver, and it serves the precision of the residual, not its
continuity. The unsymmetric
cycle as a stationary iteration is the fastest outer iteration; a Krylov method gains less than one cycle over it.

On the GPU, the patch update runs in registers, one cell per thread, and the cells of a patch exchange data by warp
shuffles. At $p=3$, a realization that stages its intermediates through shared memory is bound by a bandwidth exactly
balanced with the FP32 arithmetic it feeds; register residency avoids this bound, and the arithmetic runs close to the
FP32 peak. Fast diagonalization, merged transforms, and the reflection symmetry of the patch then remove arithmetic
until a sweep moves its data at $87\%$ of the bandwidth of a plain device copy. The solver reaches about $400$ million DoFs per second on one
A100, $2.2$ times the throughput of the mixed-precision GPU patch smoother~\cite{cui2025implementation} at
the same number of V-cycles.

The reductions of the arithmetic require Kronecker-sum patch matrices, which on Cartesian meshes means constant
coefficients and cells of equal size within a patch. Fast diagonalization is thus limited to separable problems;
for more general ones, the exact local solve can be replaced by a local $p$-multigrid~\cite{wichrowski2025local}. The smoother extends to locally refined
hierarchies by the rule that excludes refinement-edge DoFs, and the one-cell-per-thread mapping to higher degrees; both
remain to be tested, the latter against the register limit, which the $64$ coefficients of a $\mathbb{Q}_3$ cell
already approach once the patch intermediates are added.

\section*{Acknowledgments}
During the preparation of this work the author used Claude (Anthropic) to assist with drafting and editing the text.
The author reviewed and edited all content and takes full responsibility for the publication.

\bibliographystyle{siamplain}
\bibliography{literature,added_literature}

\end{document}